\documentclass{article}
\usepackage{amssymb,amsmath}
\usepackage{anysize}
\usepackage{psfrag}
\newtheorem{def.}{Definition}[section]
\newtheorem{ex}{Example}[section]

\newtheorem{prop}{Proposition}[section]

\newtheorem{lem}{Lemma}[section]

\numberwithin{table}{section}

\begin{document}
\title{Permutations Which Make Transitive Groups Primitive}
       \author{Pedro Lopes\\
        Department of Mathematics\\
        Instituto Superior T\'ecnico\\
        Technical University of Lisbon\\
        Av. Rovisco Pais\\
        1049-001 Lisbon\\
        Portugal\\
        \texttt{pelopes@math.ist.utl.pt}\\
}
\date{July 09, 2009}
\maketitle

\bigbreak

\begin{abstract} In this article we look into characterizing primitive groups in the following way. Given a primitive group we single out a subset of its generators  such that these generators alone (the so-called primitive generators) imply the group is primitive. The remaining generators ensure transitivity or comply with specific features of the group. 

We show that, other than the symmetric and alternating groups, there are infinitely many primitive groups with one primitive generator each. These primitive groups are certain Mathieu groups, certain projective general and projective special linear groups, and certain subgroups of some affine special linear groups.
\end{abstract}

\bigbreak


\bigbreak


\bigbreak

\section{Introduction} \label{sect:intro}

\noindent

We recall the basic facts and definitions concerning primitive groups in order to develop the notation. For further information we refer the reader to \cite{dm} and \cite{cameron}.

\bigbreak

Let $d$ be an integer greater than one and consider the symmetric group on $d$ letters, $\Sigma\sb{d}$, acting on the set $\Omega\sb{d} = \{ 1, 2, \dots , d \}$. We let $A\sb{d}$ stand for the alternating group. Subgroups $G$ of $\Sigma\sb{d}$ are called \emph{permutation groups}. A permutation group is said to be \emph{transitive} if for any distinct $i, i'\in \{ 1, 2, \dots , d \}$, there is a $\gamma\in G$ such that $i'=\gamma (i)$. Our notation for $G$ being a transitive subgroup of $\Sigma\sb{d}$ is
\[
G \, \,  \overset{trans}{\leq}\, \,   \Sigma\sb{d}
\]

Given a permutation group $G$, a set $B\subset \Omega\sb{d}$ is said to be a \emph{block} for $G$ if for any $\gamma\in G$, either $\gamma B=B$ or $\gamma B \cap B = \emptyset$. A transitive permutation group $G$ is called \emph{primitive} if its blocks are the trivial ones i.e., $\emptyset$, singletons, or $\Omega\sb{d}$; otherwise $G$ is called \emph{imprimitive}.

Furthermore, for a transitive permutation group, $G$, the size of a block $B$, other than the empty set, divides $d$, the number of letters being permuted; the action of $G$ on $B$ gives rise to a decomposition of $\Omega\sb{d}$ into blocks, the ``translates'' of $B$, which is called the \emph{system of blocks generated by $B$} (\cite{dm}, \cite{cameron}).

\bigbreak

In this article, we address the following situation. Given a transitive group, we look into the possibility of singling out a subset of this group such that the elements of this subset alone imply this group is primitive. Specifically, we want to find a minimal subset, $P$, of a set of generators, such that for a general candidate for a non-trivial block, $A$, the condition $\gamma A=A$ or $\gamma A\cap A = \emptyset$ breaks down for elements $\gamma\in\langle P \rangle $, the subgroup generated by $P$. We call the elements of $P$, \emph{primitive generators}, and the elements of $\langle P \rangle $, \emph{primitive permutations}.

In references \cite{neumann} or \cite{isaacs}, efforts are directed towards characterizing minimal sets of permutations which generate the symmetric or alternating groups, given that the group under study is primitive. Here we look into the possibility of constructing a set of generators for a primitive group, with a distinguished minimal subset, $P$, of primitive generators. We proceed in two steps. In the first step, the primitive generators are constructed. In the second step, generators are added to the primitive generators in order to ensure transitivity or to comply with specific features of the group.

\bigbreak

More specifically, in this article, we consider an infinite class of permutations with the following property. Each element of this class becomes a primitive generator in any transitive permutation group it may belong to. In particular, the $P$  sets here are singletons, in the notation above. Furthermore, for infinitely many primitive groups, we obtain presentations of this form i.e., with one primitive generator plus extra generators to ensure transitivity and to comply with specific features of the group. Besides the symmetric and alternating groups, these primitive groups are the Mathieu groups $M\sb{11}$, $M\sb{12}$ and $M\sb{24}$; the groups $PSL(2, p)$ and $PGL(2, p)$, for odd prime $p$; and certain subgroups of $ASL(m, 2)$, with $p+1=2\sp{m}$, where $m$ is an odd prime.

\bigbreak

Each primitive generator we will be considering here is a permutation $\alpha$ of $\Sigma\sb{d}$. The relationship between partitions of $d$ and $\alpha$ will be important below. We thus now elaborate on partitions and permutations.

\bigbreak

We henceforth use the term {\rm\bf partition of $d$} to mean a finite number of positive integers whose sum is $d$. Notation: $[n\sb{1}, n\sb{2}, \dots , n\sb{l}]$, when it is clear which integer is being partitioned. Each of the $n\sb{t}$'s is called a {\rm\bf part} of the partition $[n\sb{1}, n\sb{2}, \dots , n\sb{l}]$.

\bigbreak

Consider then a partition of $d$, $[n\sb{1}, n\sb{2}, \dots , n\sb{l}]$, and a permutation $\alpha\in\Sigma\sb{d}$ with
\[
\alpha = (i\sp{1}\sb{1},  \dots , i\sp{1}\sb{n\sb{1}})(i\sp{2}\sb{1},  \dots , i\sp{2}\sb{n\sb{2}}) \cdots (i\sp{l}\sb{1},  \dots , i\sp{l}\sb{n\sb{l}})
\]
where, for $s=1, 2, \dots , l$, the $(i\sp{s}\sb{1},  \dots , i\sp{s}\sb{n\sb{s}})$ are disjoint cycles and
\[
\Omega\sb{d} =\bigcup\sb{s=1}\sp{l} \{  i\sp{s}\sb{1},  \dots , i\sp{s}\sb{n\sb{s}} \}
\]

We will let $S\sp{s}:=\{ i\sp{s}\sb{1},  \dots , i\sp{s}\sb{n\sb{s}} \}$, for each $s=1, 2, \dots , l$. Furthermore, $[n\sb{1}, n\sb{2}, \dots , n\sb{l}]$ will be called the {\rm\bf partition of $d$ associated with $\alpha$}, with the following notation:
\[
P\sb{\alpha}(d) := [n\sb{1}, n\sb{2}, \dots , n\sb{l}]
\]

\bigbreak

{\bf In the sequel we use the symbols $d$, $\Omega\sb{d}$, $\Sigma\sb{d}$, $A\sb{d}$, $\alpha $, $[n\sb{1}, \dots , n\sb{l}]$, $(i\sp{s}\sb{1},  \dots , i\sp{s}\sb{n\sb{s}})$, $S\sp{s}$ and $P\sb{\alpha}(d)$ with the meaning we have introduced above}.

\bigbreak

For each primitive generator we will be considering here, the way the lengths of its disjoint cycles partition $d$ is crucial for our results. In particular, the following definitions will be instrumental for our results.

\bigbreak

\begin{def.}\label{def:relprimepart} A partition of $d$ into distinct integers, $[n\sb{1}, n\sb{2}, \dots , n\sb{l}]$, is called a {\rm\bf relatively prime partition of $d$} if the parts are mutually prime i.e., $(n\sb{t}, n\sb{t'})=1$ for each $t\neq t'\in\{ 1, 2, \dots , l  \}$.
\end{def.}

\bigbreak

\begin{def.}\label{def:$m$-part} A partition of $d$ into distinct integers, $[n\sb{1}, n\sb{2}, \dots , n\sb{l}]$,  is called an {\rm\bf $m$-partition of $d$} if there is a factorization $d=m k$ with integers $1<m, k<d$, along with a partition of $l$, say $[l\sb{1}, l\sb{2},  \dots , l\sb{k}]$, with integer $1<k<l$, such that, possibly after relabeling the parts,
\[
m=\sum\sb{j=1}\sp{L\sb{1}}n\sb{j}=\sum\sb{j=L\sb{1}+1}\sp{L\sb{2}}n\sb{j}= \cdots = \sum\sb{j=L\sb{k-1}+1}\sp{L\sb{k}}n\sb{j}
\]
where $L\sb{r}=\sum\sb{t=1}\sp{r}l\sb{t}$, for $r=1, 2, \dots , k$.
\end{def.}

\bigbreak

Roughly speaking, an $m$-partition of an integer $d$ is a partition such that $m$ is a non-trivial divisor of $d$ and we can gather the rearranged parts in sums adding up to $m$. Example: $d=10$, $n\sb{1}=2, n\sb{2}=5, n\sb{3}=3$; $m=5$.

\bigbreak

\begin{def.}\label{def:special$m$-part} A partition of $d$ into distinct integers, $[n\sb{1}, n\sb{2}, \dots , n\sb{l}]$,  is called a {\rm\bf special $m$-partition of $d$} if the following holds.

After relabeling, $n\sb{l}$ is the largest part and there is a common factor $m$ of $n\sb{l}$ and $d$ $(1<m<n\sb{l})$. There is, also, a partition of $l-1$, say $[l\sb{1}, l\sb{2}, \dots , l\sb{k}]$, with integer $1<k<l-1$, such that,
\[
m=\sum\sb{j=1}\sp{L\sb{1}}n\sb{j}=\sum\sb{j=L\sb{1}+1}\sp{L\sb{2}}n\sb{j}= \cdots = \sum\sb{j=L\sb{k-1}+1}\sp{L\sb{k}}n\sb{j}
\]
where $L\sb{r}=\sum\sb{t=1}\sp{r}l\sb{t}$, for $r=1, 2, \dots , k$.
\end{def.}

\bigbreak

Roughly speaking, in a special $m$-partition the largest $n\sb{t}$ is divisible by a factor, $m$, which also factors $d$. Moreover, the remaining parts can be rearranged such that we can gather the rearranged parts in sums adding up to $m$. Example: $d=15$, $n\sb{1}=2, n\sb{2}=3, n\sb{3}=10$; $m=5$. Here is a more elaborate example which is also a relatively prime partition: $d=185$, $n\sb{1}=1$, $n\sb{2}=2$, $n\sb{3}=5$, $n\sb{4}=7$, $n\sb{5}=17$, $n\sb{6}=19$, $n\sb{7}=23$, $n\sb{8}=111$; $m=37$.
\begin{align*}
185 &= 111 + ( 2 + 5 + 7 + 23 ) + ( 1 + 17 + 19 ) \\
&= 3\cdot 37 + 37 + 37
\end{align*}

\bigbreak

The main result in this work is the statement that a permutation $\alpha \in \Sigma\sb{d}$ with $P\sb{\alpha}(d)$ a relatively prime partition which is not an $m$-partition nor a special $m$-partition, constitutes the minimal set of primitive generators of a primitive group. Moreover, there is an infinite class of such groups other than the symmetric groups and the alternating groups. Specifically,

\bigbreak

{\bf Theorem}  \quad Let
\[
\alpha \in G \, \,  \overset{trans}{\leq}\, \,   \Sigma\sb{d}
\]
and let the decomposition of $\alpha$ into disjoint cycles be:
\[
\alpha = (i\sp{1}\sb{1},  \dots , i\sp{1}\sb{n\sb{1}})(i\sp{2}\sb{1},  \dots , i\sp{2}\sb{n\sb{2}}) \cdots (i\sp{l}\sb{1},  \dots , i\sp{l}\sb{n\sb{l}})
\]
with $S\sp{r}:=\{ i\sp{r}\sb{1},  \dots , i\sp{r}\sb{n\sb{r}} \}$, for each $r=1, 2, \dots , l$.

Furthermore, in the notation above,
\[
P\sb{\alpha}(d) := [n\sb{1}, n\sb{2}, \dots , n\sb{l}]
\]

If
\begin{itemize}
\item $l=2$ and $P\sb{\alpha}(d)$ is a relatively prime partition; or
\item $l\geq 3$ and $P\sb{\alpha}(d)$ is a relatively prime partition which is not an $m$-partition, nor a special $m$-partition,
\end{itemize}
then $G$ is primitive.

\bigbreak

\bigbreak

{\bf Corollary 1}

\begin{itemize}
\item For $l\geq 3$,  $G$ is either $\Sigma\sb{d}$ or $A\sb{d}$.
\item For $l=2$
\begin{itemize}
\item If $P\sb{\alpha}(d)=[n\sb{1}, n\sb{2}]$ with $n\sb{1}, n\sb{2} > 1$, then $G$ is either $\Sigma\sb{d}$ or $A\sb{d}$
\item If $P\sb{\alpha}(d) = [1, d-1]$, then $G$ may be distinct from $\Sigma\sb{d}$ or $A\sb{d}$
\end{itemize}
\end{itemize}

\bigbreak

\bigbreak

{\bf Corollary 2} \quad There are infinitely many primitive groups (other than the symmetric and alternating groups) each one with a one-element set of primitive generators. Letting $\alpha$ denote this primitive generator, $P\sb{\alpha}(d) = [1, d-1]$.

\bigbreak

\bigbreak

The Theorem and the Corollaries will be proved in Section \ref{sect:l=2} (the $l=2$ instance of the Theorem) and Section \ref{sect:proof}.

Moreover, it will be shown below that, if $P\sb{\alpha}(d)$ is a relatively prime (special) $m$-partition of $d$ associated with $\alpha (\in G)$, the transitive group $G$ at issue may be either primitive or imprimitive. On the other hand, when $l=2$, there are no relatively prime (special) $m$-partitions of $d=n\sb{1}+n\sb{2}$.

As a matter of fact, it was a special case of the $l=2$ situation, presented as an example in the PhD thesis of Natalia Bedoya (\cite{bed}), that prompted us into writing this article. In her PhD thesis, she shows that any transitive subgroup $G$ of $\Sigma\sb{d}$ containing a permutation $\alpha = (1,  2,  \dots , d-1)(d)$ (i.e., the $l=2$, $P\sb{\alpha}(d) = [1, d-1]$ situation, here) is primitive. Bedoya does so by observing that any proper subset, $A\subset\Omega\sb{d}$ containing at least two elements, one of which is $d$, is such that $\alpha A \neq A$ and $\alpha A \cap A \neq \emptyset$. Bedoya's interest in primitive groups has to do with a theorem of hers in the thesis, stating that branched coverings of surfaces, surjective at the level of the fundamental groups, are indecomposable if and only if the associated permutation groups are primitive. In this respect, it is relevant to be able to construct primitive and imprimitive groups.

\bigbreak

\bigbreak

This article is organized as follows. In Section \ref{sect:fundlemma} we state and prove a Fundamental Lemma. This allows us to settle the $l=2$ instance of the Theorem in Section \ref{sect:l=2}. In Section \ref{sect:relrpime$m$-part} we present examples illustrating the relevance of $m$-partitions and special $m$-partitions which are also relatively prime partitions. In Section \ref{sect:proof} we prove the $l\geq 3$ instance of the Theorem and the Corollaries. In Section \ref{sect:finrems} we suggest directions for further research.

\section{The fundamental lemma}\label{sect:fundlemma}

\noindent

\bigbreak

In this Section we prove the Fundamental Lemma. This lemma concerns a permutation $\alpha$ in a transitive subgroup $G$ of $\Sigma\sb{d}$. We recall that the $S\sp{r}$'s are the sets of elements moved by the disjoint cycles $\alpha$ splits into. Roughly speaking, the Fundamental Lemma states that if a proper subset $A\subset\Omega\sb{d}$ overlaps two of the $S\sp{r}$'s  whose lengths are relatively prime and without properly containing one of them, then $A$ cannot be a block.  In  subsequent Sections we will see that the combinatorics of the  partitions of $d$ under study imply that candidates for non-trivial blocks will satisfy the hypothesis of the Fundamental Lemma in a number of situations. Applying the Fundamental Lemma will then allow us to disregard these situations.

\bigbreak

\begin{lem}[Fundamental Lemma]\label{lem:fund} Let
\[
\alpha \,\,  \in \,\, G \, \,   \overset{trans}{\leq}\, \,   \Sigma\sb{d}
\]
and let the decomposition of $\alpha$ into disjoint cycles be
\[
\alpha = (i\sp{1}\sb{1}, \dots , i\sp{1}\sb{n\sb{1}})(i\sp{2}\sb{1},  \dots , i\sp{2}\sb{n\sb{2}}) \cdots (i\sp{l}\sb{1},  \dots , i\sp{l}\sb{n\sb{l}})
\]
and $S\sp{r} = \{ i\sp{r}\sb{1},  \dots , i\sp{r}\sb{n\sb{r}}  \}$, for $1\leq r \leq l$.

If $A$ is a proper subset of $\Omega\sb{d}$ such that there are distinct $s$ and $t$ in $\{ 1, \dots , l \}$ with $(n\sb{s}, n\sb{t})=1$ and
\[
\emptyset \neq \bigr( S\sp{s}\setminus A \bigl)\, \subsetneq S\sp{s} \qquad \text{ and } \qquad S\sp{t}\cap A \neq \emptyset
\]
then $A$ is not a block for $G$.
\end{lem}Proof: Let $A$ be a subset of $\Omega\sb{d}$ as in the statement. Pick $a\sb{s}\in S\sp{s}\cap A$ and $a\sb{t}\in S\sp{t}\cap A$. Then,
\[
(i\sp{t}\sb{1},  \dots , i\sp{t}\sb{n\sb{t}})\sp{n\sb{t}}(i) = i \qquad \text{ for any } i\in S\sp{t}
\]
and in particular for $i=a\sb{t}$.

Since, by hypothesis, $(n\sb{s}, n\sb{t})=1$, then
\[
(i\sp{s}\sb{1},  \dots , i\sp{s}\sb{n\sb{s}})\sp{n\sb{t}}
\]
generates the cyclic group
\[
\langle  \, (i\sp{s}\sb{1},  \dots , i\sp{s}\sb{n\sb{s}})\,   \rangle
\]
Consequently, there is an integer $n'$ such that
\[
\bigl( (i\sp{s}\sb{1},  \dots , i\sp{s}\sb{n\sb{s}})\sp{n\sb{t}}\bigr) \sp{n'} (a\sb{s}) \in S\sp{s}\setminus A
\]

In this way, the permutation $\alpha \sp{n\sb{t}n'}$ belongs to $G$ and is such that there are distinct elements $a\sb{s}$ and $a\sb{t}$ in $A$ such that
\[
\alpha\sp{n\sb{t}n'} (a\sb{s}) \notin A  \qquad \text{ and } \qquad \alpha\sp{n\sb{t}n'} (a\sb{t}) \in A
\]
Therefore, $A$ is not a block for $G$. $\hfill \blacksquare $

\bigbreak

\section{Application: the $l=2$ instance of the Theorem}\label{sect:l=2}

\noindent

We start this Section with the following observation.

\begin{lem} Let $n\sb{1}$ and $n\sb{2}$ be positive integers with $(n\sb{1}, n\sb{2})=1$. Then $(n\sb{1}, n\sb{1}+n\sb{2})=1=(n\sb{2}, n\sb{1}+n\sb{2})$.\end{lem}Proof: Let $d=n\sb{1}+n\sb{2}$. If $(n\sb{1}, d)>1$ then there would be an integer $k>1$ such that $n\sb{1}=kn\sb{1}'$  and $d=kd'$. Then
\[
kd' = d = n\sb{1}+n\sb{2} = kn\sb{1}' + n\sb{2}
\]
and thus
\[
n\sb{2} = k(d'-n\sb{1})
\]
which conflicts with $(n\sb{1}, n\sb{2})=1$.
$\hfill \blacksquare$

\bigbreak

Therefore, given any relatively prime partition $[n\sb{1}, n\sb{2}]$ of $d$, any factor of $d$ is relatively prime to the parts $n\sb{1}$ and  $n\sb{2}$. In particular, if $\Omega\sb{d}$ is partitioned by proper subsets of equal size $m \mid d$, then $m\nmid n\sb{1}$ and $m\nmid n\sb{2}$. Consequently, one of these sets will intersect both $\{ i\sp{1}\sb{1}, \dots , i\sp{1}\sb{n\sb{1}} \}$ and $\{ i\sp{2}\sb{1}, \dots , i\sp{2}\sb{n\sb{2}} \}$ without containing both of them. Hence, applying the Fundamental Lemma \ref{lem:fund}, one obtains the $l=2$ instance of the Theorem:

\begin{prop}\label{prop:l=2}Let $n\sb{1}$ and $n\sb{2}$ be positive and relatively prime integers such that $d=n\sb{1}+n\sb{2}$. If
\[
\alpha \,\,  \in \,\, G \, \,   \overset{trans}{\leq}\, \,   \Sigma\sb{d}
\]
with $P\sb{\alpha}(d) = [n\sb{1}, n\sb{2}]$, then $G$ is primitive.
\end{prop}

\bigbreak

\section{The relevance of $m$-partitions and special $m$-partitions}\label{sect:relrpime$m$-part}

\noindent

The purpose of Example \ref{ex:specialsit} below is to illustrate that a permutation $\alpha$ whose associated partition is a relatively prime partition which is also an $m$-partition, can give rise both to primitive and to imprimitive groups containing it.

\bigbreak

\begin{ex}\label{ex:specialsit}
\end{ex}

We consider a permutation $\alpha$ whose associated partition is a relatively prime partition which is also an $m$-partition:
\[
\alpha = (1, 2)(3, 4, 5)(6, 7, 8, 9, 10) \in \Sigma\sb{10}
\]
For any transitive subgroup of $\Sigma\sb{10}$ containing this permutation, the only admissible non-trivial blocks are
\[
A\sp{1} = \{ 1, 2, 3, 4, 5  \}   \qquad \text{ and } \qquad A\sp{2} = \{ 6, 7, 8, 9, 10  \}
\]
since the only non-trivial divisors of $10$ are $2$ and $5$ and since blocks of size $2$ would not cover completely the sets $\{ 3, 4, 5 \}$ or $\{ 6, 7, 8, 9, 10 \}$, therefore conflicting with the Fundamental Lemma \ref{lem:fund}. Furthermore, any other choice of sets $A\sp{1}$ and $A\sp{2}$ of size $5$ would conflict with the Fundamental Lemma.

Consider then
\[
G\sb{1}\cong\langle  \alpha , (2, 3)(5, 6)   \rangle   \, \,  \overset{trans}{\leq}\, \,   \Sigma\sb{10}
\]
This is a primitive group since the permutation $(2, 3)(5, 6)$ fixes $4$ and moves $5$ to $6$. Using \cite{gap} we find its order is $10!$. Hence, this group is $\Sigma\sb{10}$. Alternatively, using remarks $1$. and $2$. in the proof of Corollary $1$ (see below), we arrive at the same conclusion.

\bigbreak

Consider now
\[
G\sb{2}\cong\langle  \alpha , (1, 6)(2, 7)(3, 8)(4, 9)(5, 10)   \rangle   \, \,  \overset{trans}{\leq}\, \,   \Sigma\sb{10}
\]

It is easy to see that $\alpha$ and $(1, 6)(2, 7)(3, 8)(4, 9)(5, 10)$ act on the sets $ A\sp{1}$ and $A\sp{2}$ by fixing them or permuting them. Since any element in $G$ is the composition of these elements and/or their inverses, the same is true for any element of $G$. Thus, $\{ A\sp{1}, A\sp{2}\}$ is a system of non-trivial blocks for the group $G\sb{2}$. $G\sb{2}$ is then imprimitive.

\bigbreak

Analogously, the purpose of Example \ref{ex:specialsitbis} below is to illustrate that a permutation $\alpha$ whose associated partition is a relatively prime partition which is also a special $m$-partition can give rise to both primitive and imprimitive groups containing it.

\bigbreak

\begin{ex}\label{ex:specialsitbis}
\end{ex}

In the notation above,

\[
d=30; \qquad n\sb{1}=25, \quad n\sb{2}=2, \quad n\sb{3}=3
\]

\[
\alpha = (1,  2,  3,  4,  5,  6,  7,  8,  9,  10, 11, 12, 13, 14, 15, 16, 17, 18, 19, 20, 21, 22, 23, 24, 25)(26, 27)(28, 29, 30)
\]

\[
\beta = (25, 26)(27, 28)
\]

\[
\gamma = (1, 26)(6, 27)(11, 28)(16, 29)(21, 30)
\]

\bigbreak

Consider then the groups

\[
G\sb{1}\cong\langle \alpha , \beta \rangle \qquad \qquad \text{ and }\qquad \qquad  G\sb{2}\cong\langle \alpha , \gamma \rangle
\]

Clearly, both groups are transitive. The non-trivial divisors of $30$ are $2$, $3$, $5$, $6$, $10$, and $15$. Then, applying the Fundamental Lemma to $\alpha$, we see that a non-trivial block for these groups cannot be of sizes $2$, $3$, $6$, $10$, or $15$. We are left with $5$ for the size of candidates for non-trivial blocks of these groups. Again applying the Fundamental Lemma to $\alpha$, one of the blocks of size $5$ has to contain $26, 27, 28, 29, 30$. Let us, then, set
\[
A\sp{6} = \{ 26, 27, 28, 29, 30 \}
\]
Since $\beta (29) = 29\in A\sp{6}$ while $26\in A\sp{6}$ and $\beta (26) = 25\notin A\sp{6}$, then $A\sp{6}$ is not a block for $G\sb{1}$. Hence, $G\sb{1}$ is primitive. Since $|G\sb{1}| = 30!$ (using \cite{gap}) then $G\sb{1} \cong \Sigma\sb{30}$ (or via remarks $1$. and $2$. in the proof of Corollary $1$).

\bigbreak

Consider now, the group $G\sb{2}$ and, along with block $A\sp{6}$, above, the blocks

\[
A\sp{1}=\{ 1, 6, 11, 16, 21  \} \qquad \qquad A\sp{2}=\{ 2, 7, 12, 17, 22  \} \qquad \qquad A\sp{3}=\{ 3, 8, 13, 18, 23 \}
\]

\[
A\sp{4}=\{ 4, 9, 14, 19, 24 \} \qquad \qquad A\sp{5}=\{ 5, 10, 15, 20, 25 \}
\]

Then, $\alpha$ fixes $A\sp{6}$ while permutes the remaining blocks among themselves. Also, $\gamma$ maps $A\sp{1}$ and $A\sp{6}$ into each other while fixing the other blocks. Then, the $A\sp{r}$'s constitute a system of non-trivial blocks for $G\sb{2}$. Therefore, $G\sb{2}$ is imprimitive.

\bigbreak

\section{Proof of the Theorem and of the Corollaries}\label{sect:proof}

\noindent

Proof (of the Theorem): Let $A$ be a non-trivial block for $G$ which is contained in $S\sp{s\sb{0}}$, $A\subseteq S\sp{s\sb{0}}$. Then, $A$ is also a non-trivial block for the group $\langle  (\, i\sp{s\sb{0}}\sb{1}, \, \dots \, , i\sp{s\sb{0}}\sb{n\sb{s\sb{0}}})   \rangle$. Hence $|A|\, \big| n\sb{s\sb{0}}$.

Assume there is a translate of $A$ which is contained in $S\sp{t}$, for some $t\neq s\sb{0}$. With the same argument we have just used, $|A|\, \big| n\sb{t}$. Hence, $(n\sb{s\sb{0}}, n\sb{t})\neq 1$ contradicting the hypothesis.

Therefore, if $A$ is contained in $S\sp{s\sb{0}}$ then either $A=S\sp{s\sb{0}}$ or the union of $A$ with a number of translates of $A$ equals $S\sp{s\sb{0}}$. Each of the remaining translates of $A$ overlaps at least two distinct $S\sp{t}$'s (for if $A\subset S\sp{u}$ then $|A|\big| n\sb{u}$ and $(n\sb{s\sb{0}}, n\sb{u})>1$ contradicting the hypothesis).

Analogously, begin with a non-trivial block $A$ which overlaps two distinct $S\sp{t}$'s. Assume further that one of its translates is contained in one $S\sp{s\sb{0}}$. Then no other translate of $A$ is contained in an $S\sp{t}$ other than $S\sp{s\sb{0}}$, for otherwise, $(n\sb{\sb{s\sb{0}}}, n\sb{t})\neq 1$.

In this way, given a non-trivial block $A$ for $G$, some of its translates may partition one $S\sp{s\sb{0}}$ while each of the remaining translates overlaps at least two distinct $S\sp{t}$'s and does not overlap $S\sp{s\sb{0}}$.

Let us now analyze these translates which overlap two distinct $S\sp{t}$'s. Suppose there is a translate of the non-trivial block $A$, call it $A'$, with two $t\neq t'$ such that $\emptyset\neq\bigl( S\sp{t}\setminus A'\bigr) \subsetneq S\sp{t}$ and $A'\cap S\sp{t'}\neq \emptyset$. Then the Fundamental Lemma \ref{lem:fund} implies $A'$ cannot be a block therefore concluding the proof in this case.

There is, then, no translate of $A$, including $A$ itself, satisfying the hypothesis of the Fundamental Lemma.  So,
\begin{enumerate}
\item either the blocks generated by the non-trivial block $A$ are such that each is the union of a number of $S\sp{t}$'s, distinct $S\sp{t}$'s for distinct blocks;
\item or some of the blocks generated by $A$ partition an $S\sp{s\sb{0}}$, and each of the remaining blocks equal the union of distinct $S\sp{t}$'s, distinct $S\sp{t}$'s for distinct blocks. In particular, this implies that $n\sb{s\sb{0}}$ is the largest of the $n\sb{t}$'s. Otherwise, either a translate of $A$ would be contained in another $S\sp{t\sb{0}}$ or it would overlap distinct $S\sp{t}$'s without containing them. These situations have already been discarded.
\end{enumerate}
The first case implies the partition $[n\sb{1}, \dots , n\sb{l}]$ is an $m$-partition and the second case implies $[n\sb{1}, \dots , n\sb{l}]$ is a special $m$-partition. In each case there is a conflict with the hypothesis. The proof is complete.
$\hfill \blacksquare$

\bigbreak

\bigbreak

Proof (of Corollary 1) We start by listing some facts in order to prove the statement of Corollary $1$.
\begin{enumerate}
\item If $m\in P\sb{\alpha}(d)$, then there is a cycle of length $m$ in $G$. This is true because $P\sb{\alpha}(d)$ is a relatively prime partition; thus raising $\alpha$ to the products of the lengths of all its cycles except the one of length $m$, yields a cycle of length $m$.
\item If a primitive group contains a $2$-cycle it is a symmetric group; if it contains a $3$-cycle it is either an alternating group or a symmetric group (\cite{isaacs}).
\item A \emph{$k$-transitive group} is a group of permutations such that whenever two $k$-tuples of distinct letters are fixed, there is a permutation in the group that maps one $k$-tuple to the other, coordinatewise. A primitive group of permutations on $d$ letters containing an $m$-cycle with $1<m<d$ is $(d-m+1)\text{-transitive}$. This is a theorem of Marggraf; see \cite{wielandt} on page 38.
\item The $6$-  and higher transitive groups are the symmetric and alternating groups (\cite{dm}, \cite{cameron}).
\end{enumerate}
\bigbreak
In this way, if we find a part $m$ in $P\sb{\alpha}(d)$ such that $m=2$ or $3$,  or with $d-m\geq 5$ then the proof is complete. We now look into the different possibilities.

\bigbreak

If $l=2$ then either $P\sb{\alpha}(d)=[1, d-1]$ or $P\sb{\alpha}(d)=[m, d-m]$ with $m>1$. For the latter,
\begin{itemize}
\item If $m=2$ or $m=3$, ditto.
\item For $m=4$, $d$ has to be odd for otherwise the two parts of the partition would not be relatively prime. If $d=5$, this is an instance of the $P\sb{\alpha}(d)=[1, d-1]$ situation. If $d=7$, then $d-m=3$, which means the second part gives rise to a $3$-cycle. If $d\geq 9$, then $d-4\geq 5$.
\item If $m\geq 5$, then either we have an instance of the $P\sb{\alpha}(d)=[1, d-1]$ situation, or $d-m \in \{ 2, 3 , 5, 6, 7, \dots  \}$ (in which case the proof is complete), or $d-m = 4$. In the latter case, the other part of the partition is $d-m$ which satisfies $d-(d-m) = m \geq 5$.
\end{itemize}
\bigbreak
If $l\geq 3$ then the smallest parts of the partition can be $1$, $2$ and $3$ in which case the proof ends.
\begin{itemize}
\item If there is an $m=4$ then the other parts may be as small as $1$ and $3$; or one of them is $1$ and the others greater than $4$, or all the other parts are greater than $4$ which implies that $d$ minus one of these is greater than or equal to $5$.
\item If one of the parts is $m\geq 5$ then the smaller parts include $2$, $3$ or $m'$ with $3 < m' < m$. If $d-m =4$, then $d - m' \geq 5$. The remaining situations are dealt with as above. This completes the proof.
\end{itemize}
$\hfill \blacksquare$

\bigbreak

We remark that the symmetric and alternating groups can be realized in this way i.e., with a primitive generator plus a second generator ensuring transitivity. For even $d\geq 4$, take $\alpha = (1)(2, 3, \dots , d)$ and $\beta = (1, 2)(3, 4)$ for the alternating group, $A\sb{d}$, and $\beta = (1, 2)$ for the symmetric group, $\Sigma\sb{d}$.

For odd $5\leq d = 2k+1$, $k>3$, take $\alpha = (1, 2, \dots , k-1)(k, k+1, \dots , 2k)(2k+1)$.  $A\sb{3}$, $A\sb{5}$ and $A\sb{7}$ can be realized in this way provided we take for the primitive generator $\alpha = (1, 2, 3)$, $\alpha =(1, 2, 3, 4, 5)$, $\alpha = (1, 2, 3, 4, 5, 6, 7)$, respectively. The same primitive generators apply for $\Sigma\sb{3}$, $\Sigma\sb{5}$, and $\Sigma\sb{7}$. We leave the details to the reader.

\bigbreak

Proof (of Corollary 2) In \cite{zieschang} we find the list of primitive groups of permutations on $p+1$ letters and containing a $p$-cycle, for prime $p$ (other than the symmetric and alternating groups). These are
\begin{itemize}
\item The Mathieu groups $M\sb{11}$, $M\sb{12}$ and $M\sb{24}$;
\item The subgroups $G$ of $ASL(m, 2)$, where $G$ is a semidirect product of the translation group $T(2, m)$ with a group $U$, with $\mathbb{Z}\sb{p}\leq U \leq SL(m, 2)$ (where $p+1=2\sp{m}$);
\item The groups $PSL(2, p)$ and the groups $PGL(2, p)$, for odd prime $p$.
\end{itemize}
In particular, note that, for prime $p\geq 7$, the groups $PSL(2, p)$ are all distinct from the symmetric and alternating groups (\cite{artin}). This completes the proof.
$\hfill \blacksquare$

\subsection{Presentations of some of the groups in Corollary $2$ with the primitive generator $\alpha$}\label{subsect:cor2}

\noindent

In this Subsection we provide presentations of some of the primitive groups in Corollary $2$ with a primitive generator $\alpha$ with $P\sb{\alpha}(d) = [d-1, 1]$.

\bigbreak

\bigbreak

The Mathieu group $M\sb{12}$. The generators are (adapted from \cite{conder}):

\begin{align*}
\alpha &= (1, 2, 3, 5, 6, 8, 9, 11, 10, 7, 4), \\
\beta &=  (3, 4)(6, 7)(9, 10)(11, 12)
\end{align*}

\bigbreak

\bigbreak

The Mathieu group $M\sb{24}$. The generators are (\cite{conway}):

\begin{align*}
\alpha &= (1, 2, 3, 4, 5, 6, 7, 8, 9, 10, 11, 12, 13, 14, 15, 16, 17, 18, 19, 20, 21, 22, 23), \\
\beta &= (16, 8, 15, 6, 11, 21, 18, 12, 23, 22, 20)(4, 7, 13, 2, 3, 5, 9, 17, 10, 19, 14), \\
\gamma &= (24, 1)(16, 4)(8, 14)(15, 19)(6, 10)(11, 17)(21, 9)(18, 5)(12, 3)(23, 2)(22, 13)(20, 7), \\
\delta &= (15, 18, 12, 20, 23)(21, 11, 8, 6, 22)(19, 5, 3, 7, 2)(9, 17, 14, 10, 13)
\end{align*}

\bigbreak

\bigbreak

The group $PSL(2, 7)$. The generators are (\cite{alperinbell}):

\begin{align*}
\alpha &= (1, 2, 3, 4, 5, 6, 7), \\
\beta &= (2, 5, 3)(4, 6, 7), \\
\gamma &= (8, 1)(2, 4)(3, 6)(5, 7)
\end{align*}

\bigbreak

\section{Final remarks}\label{sect:finrems}

\noindent

\bigbreak

The present work is an attempt at conceiving primitive groups as being yielded by sets of generators with a particular subset consisting of the so-called primitive generators i.e., generators such that they alone imply the group is primitive. The remaining generators ensure transitivity or comply with specific features of the primitive group under study.

In our examples there is only one primitive generator. This is a permutation whose associated partition is a relatively prime partition, which is neither an $m$-partition nor a special $m$-partition. We succeed in characterizing in this way, primitive groups other than the symmetric or the alternating groups. These are the Mathieu groups $M\sb{11}$, $M\sb{12}$ and $M\sb{24}$; the groups $PSL(2, p)$ and $PGL(2, p)$, for odd prime $p$; and certain subgroups of $ASL(m, 2)$, with $p+1=2\sp{m}$, where $m$ is an odd prime.

In this way, here are some questions we would like to answer.

Is it always possible to identify the primitive generators of a primitive group? Do we need only one primitive generator per primitive group or certain primitive groups require more than one? Can we characterize the class of permutations which are potential primitive generators?

Moreover, we would like to study further the relatively prime partitions which are $m$-partitions or special $m$-partitions. Also, we would like to consider other sorts of partitions like those whose parts come with multiplicity greater than one but the distinct numbers the parts take on are all mutually prime.

\bigbreak

In a different direction, we would like to study the interplay between primitive groups and indecomposability of branched coverings of surfaces as in Bedoya's \cite{bed}. In particular, what are the primitive groups here and how can we identify them.

\bigbreak

We plan to address these and other issues in future work.

\bigbreak

\subsection{Acknowledgements}

\noindent

The author acknowledges support by {\em Programa Operacional
``Ci\^{e}ncia, Tecnologia, Inova\c{c}\~{a}o''} (POCTI) of the {\em
Funda\c{c}\~{a}o para a Ci\^{e}ncia e a Tecnologia} (FCT)
cofinanced by the European Community fund FEDER.  He also thanks
the staff at IMPA and especially his host, Marcelo Viana, for
hospitality during his stay at this Institution.
The package \cite{gap} was used in the calculations which led to this article.

\bigbreak



\begin{thebibliography}{99}


\bibitem{alperinbell}
        J. L. Alperin, R. B. Bell, \emph{Groups and Representations},
        Graduate Texts in Mathematics, {\bf162},
        Springer Verlag,  1995


\bibitem{artin}
        E. Artin, \emph{The orders of the linear groups},
        Comm. Pure Appl. Math. {\bf8} (1955), 355--365




\bibitem{bed}
        N. Bedoya, \emph{Revestimentos ramificados e o problema da decomponibilidade}, (Branched coverings and their decomposability), PhD thesis, Universidade de S\~ao Paulo, S\~ao Paulo, SP, Brazil, June, 2008


\bibitem{cameron}
        P. J. Cameron, \emph{Permutation groups}, London Mathematical Society Student Texts, {\bf45}, Cambridge University Press, Cambridge, 1999



\bibitem{conder}
        M. Conder, \emph{Generating the Mathieu groups and associated Steiner systems},
        Discrete Math. {\bf112} (1993), no. 1-3, 41--47



\bibitem{conway}
        J. H. Conway, N. J. A. Sloane, \emph{Sphere packings, lattices and groups}, 3rd edition,  Grundlehren der Mathematischen Wissenschaften, {\bf290}, Springer-Verlag, New York, 1999


\bibitem{dm}
        J. D. Dixon, B. Mortimer, \emph{Permutation groups},
        Graduate Texts in Mathematics, {\bf163},
        Springer Verlag,  1996



\bibitem{gap}
        The GAP Group, GAP --- Groups, Algorithms, and Programming,
        Version 4.4.10; 2007
        (http://www.gap-system.org)



\bibitem{isaacs}
        I. M. Isaacs, T. Zieschang, \emph{Generating symmetric groups},
        Am. Math. Monthly {\bf 102} (1995), no. 8, 734-739


\bibitem{neumann}
        P. M. Neumann, \emph{Primitive permutation groups containing a cycle of prime-power length},
        Bull. London. Math. Soc. {\bf 7} (1975), 298-299


\bibitem{wielandt}
        H. Wielandt, \emph{Finite permutation groups}, Academic Press, New York-London, 1964


\bibitem{zieschang}
        T. Zieschang, \emph{Primitive permutation groups containing a $p$-cycle},
        Arch. Math., {\bf 64} (1995), 471-474





\end{thebibliography}
\end{document}